\documentclass[12pt]{amsart}

\usepackage{amsmath}
\theoremstyle{plain}
\newtheorem{Thm}{Theorem}

\begin{document}

\title[Eigenvalue estimates]
{eigenvalue estimates and L1 energy on closed manifolds}

\author{Li Ma}

\address{Department of mathematical sciences \\
Tsinghua university \\
Beijing 100084 \\
China} \email{nuslma@gmail.com}
\dedicatory{}
\date{May 26th, 2009}

\begin{abstract}

In this paper, we study Lichnerowicz type estimate for eigenvalues
of drifting Laplacian operator and L1 and L2 energy for drifting
heat equation on closed manifolds with weighted measure. In some
sense, this study is about the eigenvalue estimate on Ricci
solitons.

{ \textbf{Mathematics Subject Classification} (2000): 35J60,
53C21, 58J05}

{ \textbf{Keywords}:  Lichnerowicz type,  eigenvalue estimate,
Kato inequality, stability}
\end{abstract}

\thanks{$^*$ The research is partially supported by the National Natural Science
Foundation of China 10631020 and SRFDP 20060003002. }
 \maketitle

\section{Introduction}
In this paper we consider two kinds of problems. One is the
Lichnerowicz type eigenvalue estimate for drifting Laplacian
operator on closed Riemannian manifold $(M,g)$ with weighted
measure $e^{-f}dv_g$, where $f$ is a given smooth function on $M$.
In some sense, this problem is about the eigenvalue estimate on
Ricci solitons. This part is a continuation of previous studies in
\cite{ML1} and \cite{ML2}. See also \cite{X} for related. The
other is about the L1 energy growth estimate for heat equation
with drifting. We shall also study the growth of the Dirichlet
energy, which can be considered as a natural extension of the
eigenvalue estimate of the drifting Laplacian operator. Our
results are Theorem \ref{thm1} and Theorem \ref{thm2} below. We
use the energy method, inspired by Lichnerowicz's estimates
\cite{A98} and Qian's work \cite{Q}. Our argument is unlike the
gradient estimate through the maximum principle trick \cite{BH97}
for Harnack quantity (see \cite{BLN}, \cite{LY},\cite{SY},
\cite{L},\cite{ML1}, \cite{ML2}, and \cite{LM}). The key
ingredient is again the Bochner type identity related to the
drifting Laplacian operator and the Bakry-Emery-Ricci tensor (see
\cite{BE} and \cite{BQ}). We remark that in the recent work of
\cite{Lim} (see Theorem 3 there), Limoncu studies another kind of
Lichnerowicz estimate for Laplacian operator.

\section{Lichnerowicz type estimate}
Let us start from the proof of standard Lichnerowicz theorem
through Bochner method (one may see \cite{A98} from the point of
the Ricci formula). Let $L=\Delta $ on the compact Riemannian
manifold of dimension $n$. Assume that $f$ is a smooth function on
$M$. Recall the following Bochner formula (see
 \cite{SY})
 $$
\frac{1}{2}L|\nabla f|^2=|D^2f|^2+(\nabla f,\nabla
Lf)+(Ric)(\nabla f,\nabla f).
 $$

Let $\lambda>0$ be an eigenvalue of $-L$. Then we have a smooth
function $u$ such that
$$
-Lu=\lambda u
$$
with $\int u^2=1$. Then $\lambda=\int |\nabla u|^2$. Assume that
$Ric\geq (n-1)k$ for some $k>0$. Then using $$|D^2u|^2\geq
\frac{1}{n}|\Delta u|^2=\frac{\lambda^2}{n}u^2$$ and
$$
(\nabla u,\nabla Lu)=-\lambda |\nabla u|^2
$$
 we have by integration on $M$,
$$
0\geq \frac{\lambda^2}{n}-\lambda ^2+(n-1)k\lambda.
$$
Hence, we have
$$
\lambda\geq nk.
$$
This is the famous Lichnerowicz Theorem for eigenvalue estimate
\cite{A98}.

The above argument can be used to study the eigenvalue problem of
the elliptic operator with drifting
$$
L_f=\Delta -\nabla f\nabla
$$
associated with the weighted volume form $e^{-f}dv$. Assume that
\begin{equation}\label{eq1}
-L_f u=\lambda u.
\end{equation}
Then again $\int u^2=1$ and
$$
\lambda =\int |\nabla u|^2.
$$
Here the integral is about the weighted volume form.

Recall the following Bochner formula (see Lemma 2.1 in
 \cite{GN})
 $$
\frac{1}{2}L_f|\nabla u|^2=|D^2u|^2+(\nabla u,\nabla
L_fu)+(Ric+D^2f)(\nabla u,\nabla u).
 $$
 Assume that
 \begin{equation}\label{ass1}
Ric+D^2f\geq \frac{|Df|^2}{nz}+A
 \end{equation}
 for some $A>0$ and $z>0$.

Note that $$(a+b)^2\geq \frac{a^2}{z+1}-\frac{b^2}{z}$$ for any
$z>0$. So, we have
$$
(\Delta u)^2=(\lambda u+\nabla f\nabla u)^2\geq
\frac{\lambda^2u^2}{z+1}-\frac{|\nabla f\nabla u|^2}{z}.
$$
Then
$$
0\geq \frac{1}{n}\int(\frac{\lambda^2u^2}{z+1}-\frac{|\nabla
f|^2|\nabla u|^2}{z})+\lambda\int |\nabla u|^2+\int
(Ric+D^2f)(\nabla u,\nabla u)
$$
and
$$
0\geq \frac{1}{n}\int(\frac{\lambda^2u^2}{z+1})-\lambda\int
|\nabla u|^2+A\int |\nabla u|^2.
$$
Hence,
$$
0\geq \frac{\lambda^2}{n(z+1)}-\lambda^2+A\lambda
$$
and
$$
\lambda \geq \frac{n(z+1)A}{(n(z+1)-1)}.
$$
In conclusion, we get
\begin{Thm}\label{thm1} Let $(M,g)$ be a compact Riemannian manifold of
dimension $n$ with assumption (\ref{ass1}) for some $z>0$ and
$A>0$. Let $u$ be a smooth solution to (\ref{eq1}) with
$\lambda>0$. Then we have,
$$
\lambda \geq \frac{n(z+1)A}{(n(z+1)-1)}.
$$
\end{Thm}
This result is the Lichnerowicz type eigenvalue estimate for the
drifting Laplacian operator.

\section{new energy bound derived from Bochner formula}

We now study the following heat equation with drifting term
\begin{equation}\label{eq0}
\partial_tu=\Delta u+B\cdot \nabla u+cu
\end{equation}
on $M\times (0,T)$. Here $B$ is a smooth vector field on $M$ and
we shall let $B=\nabla f$ later, and $c\in C^2(M)$.

We define a tensor $\nabla B$ as in \cite{GN} that
$$
\nabla B(X,Y)=(\nabla_X B,Y).
$$
We assume that there is a constant $-K$ such that
\begin{equation}\label{ass}
Ric-\nabla B\geq -K
\end{equation}
on $M$.

Let $$ L^b=\Delta+ B\cdot \nabla.
$$
 Then we have the following Bochner formula (see Lemma 2.1 in
 \cite{GN})
 $$
\frac{1}{2}L^b|\nabla f|^2=|D^2f|^2+(\nabla f,\nabla
L^bf)+(Ric-\nabla B)(\nabla f,\nabla f).
 $$

From this, we know that
$$
(\partial_t-L^b)|\nabla u|^2=-2|D^2u|^2-2(\nabla u,\nabla
(cu))-2(Ric-\nabla B)(\nabla u,\nabla u).
$$

From now on, we take $c=constant$. Note that by Kato's inequality,
\begin{equation}\label{kato}
|D^2u|^2=|\nabla du|^2\geq |\nabla |du||^2\geq \frac{|\nabla
|\nabla u|^2|^2}{4(|\nabla u|^2+\epsilon)}.
\end{equation}

Let $\phi(x)$ be a monotone increasing function, i.e, $\phi'>0$,
and let $H(x)=\log \phi'(x)$. It is obvious that
$H'\phi'=\phi^{''}$. Let $w=\phi(|\nabla u|^2)$. Note that
$$
(\partial_t-L^b)w=\phi'(\partial_t-L^b)|\nabla
u|^2-\phi^{''}|\nabla |\nabla u|^2|^2.
$$
Then we have
$$
(\partial_t-L^b)w=\phi'[-2|D^2u|^2-2c|\nabla u|^2-2(Ric-\nabla
B)(\nabla u,\nabla u)-H'|\nabla |\nabla u|^2|^2].
$$
Hence, by (\ref{kato}), we have
$$
(\partial_t-L^b)w\leq \phi'[-(\frac{1}{2(|\nabla
u|^2+\epsilon)}+H')|\nabla |\nabla u|^2|^2+2(K-c)|\nabla u|^2].
$$
Assume that $B=-\nabla f$. Integrating by part  on $M$ with
respect to the weighted volume form, we have
$$
\partial_t\int_M w\leq \int_M \phi'[-(\frac{1}{2(|\nabla
u|^2+\epsilon)}+H')|\nabla |\nabla u|^2|^2+2(K-c)|\nabla u|^2].
$$
We now choose $\phi$ such that
$$
\frac{1}{2x}+H'(x)\geq 0.
$$
For example, we take $ \phi(x)=x^{p/2}$ for $x>0$ with $p\geq 1$
and $\phi'(x)=\frac{p}{2}x^{p/2-1}$. Then we have, by sending
$\epsilon \to 0$, that
$$
\partial_t\int_M w \leq 2(K-c)\int_M \phi'|\nabla u|^2.
$$
Let $p=1$, and we then have
$$
\partial_t\int_M |\nabla u| \leq (K-c)\int_M |\nabla u|,
$$
which gives us the $L^1$ energy bound of $|\nabla u|$. Take $p=2$
and we have
$$
\int_M |\nabla u|^2(t) \leq exp[2(K-c)t]\int_M |\nabla u|^2(0).
$$

In conclusion we have
\begin{Thm}\label{thm2} Let $(M,g)$ be a compact Riemannian manifold of
dimension $n$ with (\ref{ass}). Assume that $f$ is a smooth
function in $M$. Let $u$ be a smooth solution to (\ref{eq0}) with
$B=-\nabla f$ and $c=constant$. Then we have for $t>0$,
$$
\int_M |\nabla u|(t) \leq exp[(K-c)t]\int_M |\nabla u|(0).
$$
Here the integration on $M$ is with respect to the weighted volume
form. We also have
$$
\int_M |\nabla u|^2(t) \leq exp[2(K-c)t]\int_M |\nabla u|^2(0).
$$
\end{Thm}

We remark that in general, for $ \phi(x)=x^{p/2}$ we have the
$L^p$ control in the following way
$$
\partial_t\int_M |\nabla u|^p \leq (K-c)p\int_M |\nabla u|^p
$$
and then
$$
\int_M |\nabla u|^p(t) \leq exp[p(K-c)t]\int_M |\nabla u|^p(0).
$$

\end{document}